\def \version {February 14, 2025}

\documentclass[12pt]{article}

\usepackage{amsmath}
\usepackage{amssymb}
\usepackage{amsthm}
\usepackage{stmaryrd}
\usepackage{graphicx} 
\usepackage{tikz}
\usepackage{lineno}
\usepackage{comment}

\newcommand{\nev}[1]{{\bf\itshape (#1)}\, }

\def \ssk {\smallskip}
\def \msk {\medskip}
\def \bsk {\bigskip}
\def \nin {\noindent}

\newcommand{\floor}[1]{\lfloor #1 \rfloor }

\newtheorem{Theorem}{Theorem}
\def \btm {\begin{Theorem}}
\def \etm {\end{Theorem}}
\newtheorem{Lemma}[Theorem]{Lemma}
\def \blm {\begin{Lemma}}
\def \elm {\end{Lemma}}
\newtheorem{Problem}[Theorem]{Problem}
\def \bpm {\begin{Problem}}
\def \epm {\end{Problem}}
\newtheorem{Proposition}[Theorem]{Proposition}
\def \bpn {\begin{Proposition}}
\def \epn {\end{Proposition}}
\newtheorem{Corollary}[Theorem]{Corollary}
\def \bcr {\begin{Corollary}}
\def \ecr {\end{Corollary}}
\newtheorem{con}[Theorem]{Conjecture}
\def \bcj {\begin{con}}
\def \ecj {\end{con}}
\newtheorem{Homework}[Theorem]{Homework}
\def \bhw {\begin{Homework}}
\def \ehw {\end{Homework}}
\newtheorem{Example}[Theorem]{Example}
\def \bex {\begin{Example}\rm }
\def \eex {\end{Example}}
\newtheorem{Definition}[Theorem]{Definition}
\def \bdf {\begin{Definition}\rm }
\def \edf {\end{Definition}}
\newtheorem{Remark}[Theorem]{Remark}
\def \brm {\begin{Remark}\rm }
\def \erm {\end{Remark}}
\def \bpf {\begin{proof}}
\def \epf {\end{proof}}
\def \hsp {\hspace{0.5em}}

\def \Boxx {\,\Box\,}

\def \zzz {\mathbb{Z}}
\def \cF {\mathcal{F}}
\def \tw {\mathrm{tw}}
\def \rtw {\mathrm{rtw}}
\def \Mad {\mathrm{Mad}}

\def \xod {\chi_\mathrm{o}}
\def \soc {\chi_\mathrm{so}}
\def \sid {\alpha_\mathrm{od}}

\def \es {\varnothing}
\def \smin {\,\diagdown\,}
\def \soi {odd independent}
\def \sois {\soi\ set}
\def \soin  {odd independence number}

\def \sokk  {strong odd coloring}
\def \NP {{\sf NP}}

\begin{document}


\title{The odd independence number of graphs, I: Foundations and classical classes}
\author{Yair Caro\thanks{\hsp Department of Mathematics,
  University of Haifa-Oranim, Tivon 36006, Israel}\,,
 Mirko Petru\v sevski\thanks{\hsp Faculty of Mechanical Engineering, Ss. Cyril and Methodius University in Skopje, Macedonia}\,,
  Riste \v Skrekovski\thanks{\hsp Faculty of Mathematics and Physics, University of Ljubljana;  Faculty of Information Studies in Novo Mesto;  Rudolfovo - Science and Technology Centre Novo Mesto, FAMNIT, University of Primorska, Slovenia}\,,
   Zsolt Tuza\thanks{\hsp Alfr\'ed R\'enyi Institute of Mathematics,
    H-1053 Budapest, Re\'altanoda u.~13--15, Hungary;
     and
    Department of Computer Science and Systems Technology,
    University of Pannonia, 8200 Veszpr\'em, Egyetem u.~10, Hungary}}
\date{\small Latest update on \version}
\maketitle

\begin{abstract}
An odd independent set $S$ in a graph $G=(V,E)$ is an independent set of
 vertices such that, for every vertex $v \in V \smin S$, either
 $N(v) \cap S = \es$  or  $|N(v) \cap S| \equiv 1$ (mod 2),
 where $N(v)$ stands for the open neighborhood of $v$.
The largest cardinality of odd independent sets of a graph $G$,
  denoted  $\sid(G)$, is called the \soin\ of~$G$.

This new parameter is a natural companion to the recently introduced
 strong odd chromatic number.
A proper vertex coloring of a graph $G$ is a strong odd coloring
 if, for every vertex $v \in V(G)$, each color used in the neighborhood
 of $v$ appears an odd number of times in $N(v)$.
The minimum number of colors in a strong odd coloring of $G$
 is denoted by $\soc(G)$.

A simple relation involving these two parameters and the order $|G|$
 of $G$ is $\sid(G)\cdot\soc(G) \geq |G|$, parallel to the same on
 chromatic number and independence number.  

We develop several basic inequalities concerning $\sid(G)$,  and use
 already existing results on  strong odd coloring,  to derive
   lower bounds for odd independence in many families of graphs.

We prove that $\sid(G) = \alpha(G^2)$ holds for all claw-free graphs $G$, and apply this result to prove that determining $\sid(G)$ is in general \NP-hard (and also when restricted to line graphs).
We also present many results, using various techniques, concerning the odd independence number of cycles, paths, Moore graphs, Kneser graphs,
the complete subdivision $S(K_n)$ of $K_n$, the half graphs $H_{n,n}$, and $K_p \Boxx K_q$.
Further, we consider the odd independence number of the hypercube $Q_d$ and also of the complements of triangle-free graphs. 
Many open problems for future research are stated.
Further related results can be found in part II of this work [The odd independence number of graphs, II: Finite and infinite grids and chessboard graphs, arXiv: 2510.01897].

\bsk

\nin

\bsk

\nin
(primary),
05C76
(secondary)

\end{abstract}

\section{Introduction}

An \emph{odd independent set} $S$ in a graph $G = (V, E)$ is an independent set of vertices such that, for every vertex $v \in V \smin S$, either $N (v) \cap S = \es$ or $|N (v) \cap S| \equiv 1$ (mod 2), where $N (v)$ stands for the open neighbourhood  of $v$.
The largest cardinality of odd independent sets in $G$, denoted $\sid(G)$, is called the \emph{odd independence number} of $G$.
This new parameter is a natural companion to the recently introduced strong odd chromatic number \cite{CPST,GK-etal,KP-a24,PMF-DM26,P-a25}.
A proper vertex coloring of a graph $G$ is a strong odd coloring if, for every vertex $v \in V (G)$, each
color used in the neighborhood of $v$ appears an odd number of times in $N (v)$. The minimum number of colors in a strong odd coloring of $G$ is denoted by $\soc(G)$.
This paper is aimed at introducing and starting a systematic study of the odd independence number;
 further results on $\sid$ are given in the follow-up manuscript
  \cite{part-2}.

In Section \ref{s:basic} we present basic inequalities relating $\sid(G)$, $\soc(G)$ and other graphs parameters like $|G|$ (the order of $G$) and $\Delta(G)$ (the maximum degree of $G$), similar to such classical inequalities involving $\alpha(G)$ and $\chi(G)$.
We also supply upper bounds for $\sid(G)$ for $d$-regular graphs when $d$ is even, or when $G$ is a $K_{1,r}$-free graph, or when any two adjacent vertices in $G$ have at least $\lambda$ vertices in common.
We also mention that $\sid(G) = \alpha(G)$ holds for every
 regular bipartite graph $G$ with odd degree.
Moreover, concerning the class of claw-free graphs $G$,
 complementing the result $\soc(G) = \chi(G^2)$ from \cite{CPST}
 we prove that also $\sid(G) = \alpha(G^2)$ holds for them,
 and apply this identity to prove that
  determining $\sid(G)$ is \NP-hard, even when restricted to
  many types of line graphs.
Then we present various lower bounds concerning $\sid(G)$ using known results on $\soc(G)$ and the inequality $\sid(G) \geq |G|/\soc(G)$, which serves as a benchmark for further results.

In Section \ref{s:spec-grs} we develop estimates concerning $\sid(G)$ for many special classes of graphs, including the Petersen graph, Hoffman--Singleton graph and Kneser graphs, helping to assimilate the notion of odd independence.

In Section \ref{s:recurs} we present methods to recursively construct graphs from given ones, by means of graph operations inspired by the Cartesian product of graphs, which allow us to conclude a lower bound on the odd independence number
of the produced graph in terms of the building blocks.

In Section \ref{s:grid} we apply the constructions introduced in Section \ref{s:recurs} to compute a lower bound for $\sid( Q_d)$,  the $d$-dimensional cube, when $d$ is even.
The case of $d$ odd is trivially covered by the simple fact that for every $d$-regular bipartite graph $G$ with $d$ odd, $\sid(G) = \alpha(G) = |G|/2$.
In contrast, the case when $d$ is even is challenging.

In Section \ref{s:co-triangfree} we consider complements $\overline{G}$ of triangle-free graphs $G$.
These complementary graphs have $\alpha(\overline{G}) \leq 2$, hence also $\sid(\overline{G}) \leq 2$.
We analyse the possible values of $\sid(G)$ with respect to the diameter of $G$ and $\overline{G}$.

In Section \ref{s:maxdeg} we improve upon one of the basic inequalities presented in Section \ref{s:basic}, namely $\sid(G)  \geq |G|/(\Delta^2(G) +1)$.
Using our results concerning Moore graphs, we prove that the stronger inequality $\sid(G) \geq |G|/(\Delta^2(G) -1)$ is valid for all $\Delta \geq 3$.

In Section \ref{s:concl} we list many open problems for future research, covering structural problems as well as computational complexity problems.

\subsection{Notation}

Let $G$ be any graph (simple, undirected).
We denote its number $|V(G)|$ of vertices by $|G|$ (the order of $G$),
 and its number of edges by $e(G)$.
As usual, $\alpha(G)$ and $\chi(G)$ denote the independence number and the chromatic number of $G$, respectively.
Further standard notations are $N(v)$ (open neighborhood), $N[v]$ (closed neighborhood), and $\deg_G(v)$ or $\deg(v)$ (degree) of a vertex $v\in V(G)$, and similarly $N(X)$ and $N[X]$ for the corresponding neighborhoods (open/closed) of a set $X\subset V(G)$.
The average degree of graph $G$ is $d(G) := \frac{1}{|G|} \sum_{v\in V(G)} \deg(v) = \frac{2e(G)}{|G|}$.
In the square $G^2$ of $G$, two vertices are adjacent if and
   only if they are at distance at most 2 apart in $G$ (i.e., they are
   adjacent or have a common neighbor).

Well-known basic types of graphs are the complete graph $K_n$, the empty (edgeless) graph $E_n$, and the path $P_n$ of order $n$.
The graph operations $G\cup H$ (vertex-disjoint union), $G+H$ (complete join), $tG$ (the vertex-disjoint union of $t$ copies of $G$) are also standard.

Beside $\sid(G)$ and $\soc(G)$ we shall occasionally mention the further parameter $\xod(G)$, termed odd chromatic number, which is the smallest number $k$ of
   colors admitting a proper vertex $k$-coloring such that, for each non-isolated
   vertex $v\in V(G)$, some color occurs an odd number of times in $N(v)$.

\section{Some basic facts}
\label{s:basic}

In the next proposition we collect the basic relations
 between $\sid(G)$ and $\soc(G)$, parallel to the classical relation between $\alpha(G)$ and $\chi(G)$.

\bpn
\label{p:sid-LB}
Let $G=(V,E)$ be a graph on $n$ vertices.
\begin{itemize}
 \item[$(1)$] $\xod(G) \leq  \soc(G) \leq  \chi(G^2) \leq  \Delta(G^2) +1
   \leq  (\Delta (G))^2 +1$.
 \item[$(2)$] $\soc (G) = \chi(G^2)$ whenever $G$ is claw-free.
 \item[$(3)$] $\sid(G)  \geq  n/\soc (G) \geq  n/\chi(G^2) \geq
   n/(\Delta (G^2) +1)  \geq  n/((\Delta (G))^2 +1)$.
 \item[$(4.1)$] $( n +1)^2\!/4 \geq \sid(G) \soc(G) \geq n$,
   and both bounds are best possible.
 \item[$(4.2)$] $n +1 \geq \sid(G) + \soc(G) \geq 2\sqrt{n}$,
   and both bounds are best possible.
 \item[$(5)$] $\sid(G) \geq  \alpha (G^2)$.
 \item[$(6)$] $\sid(G)  \geq   \sum_{v\in V} 1/(\deg_{G^2} (v) +1 ) \geq
   n/ (d(G^2) +1)$ 
   
 \qquad\quad \ $\geq  n/ (d(G)\Delta (G) +1) \geq  n/( (\Delta(G))^2+1)$.
\end{itemize}
\epn
\bpf

\

\msk

\nin
(1)\quad
The only inequality that needs a proof is $\soc(G) \leq \chi(G^2)$.
Let us consider a coloring of $G^2$ with $\chi(G^2)$ colors.
Then in $G$ every vertex $v$ has all the vertices of $N(v)$ getting distinct colors (as $N[v]$ forms a clique in $G^2$), hence any color present in $N(v)$ appears once.

\msk

\nin
(2)\quad
This was proved already in the paper on strong odd coloring, \cite[Proposition 3]{CPST}.

\msk

\nin
(3)\quad
If $S_1\cup\cdots\cup S_k = V(G)$ is a partition into $k:=\soc(G)$ \sois s,
 then $|G| = |S_1|+\dots+|S_k| \leq k\cdot \sid(G) = \sid(G)\cdot \soc(G)$.
This implies the leftmost inequality.
The rest of the chain of inequalities follows from item~(1).

\msk

\nin
(4.1)--(4.2)\quad
Clearly from (3) we have $\sid(G) \soc(G) \geq n$. By the inequality
 between arithmetic and geometric means, denoting
$\sid(G)= x$ and $\soc(G)  = y$, we obtain $x +y \geq  2\sqrt{xy}
 \geq 2\sqrt{n}$, proving the right side of (4.2). 
Now choose a maximum odd independent set $S$, color all its vertices
 with the same color, and the remaining $n - |S|$ vertices with
 distinct colors, altogether  using $n - \sid(G) +1$ colors. 
This is easily seen to be a strong odd coloring. Hence $n +1 \geq 
 \sid(G) + \soc(G)$. Lastly, again by the inequality between arithmetic
  and geometric means, with the same notation as above,
  $( n +1)/2 \geq (x +y)/2 \geq \sqrt{xy}$ holds, hence
  we conclude that $( n +1)^2\!/4 \geq \sid(G) \soc(G)$.

Tightness of the lower bounds is shown by the graph $G= tK_r$, $n = tr$. 
Then $\soc(G) = r$ and $\sid(G) = t$, hence $n = \sid(G) \soc(G)$.
In particular, for $t = r$ and $n = r^2$ we obtain $\sid(G) + \soc(G)
 = 2r = 2\sqrt{n}$.

Tightness of the upper bounds is shown by the graph $G=K_r \cup tK_1$,
 $n=r+t$.
Then $\soc(G) = r$ and $\sid(G) = t+1$, hence $n+1 = \sid(G) + \soc(G)$.  
In particular, if $n=2k+1$ is odd, we take $r = k+1 = (n+1)/2$ and
 $t=k$, hence obtain $\sid(G) \soc(G) = r\cdot(t+1) = \frac{(n+1)^2}{4}$.
If $n=2k$ is even, we can take $r=t=k$ and obtain $\sid(G) \soc(G) = k^2+k
 = (k+1/2)^2 - 1/4 = \floor{\frac{(n+1)^2}{4}}$.

\msk

\nin
(5)\quad
Let $S$ be a maximum independent set in $G^2$ and consider $S$ as an
 independent  set in $G$ (clearly it is).
Then a vertex in $V \smin S$  is adjacent to at most one  vertex
 in $S$, since if it was adjacent to two vertices $x ,y \in S$ then
 $x , y$ are not independent in $G^2$ and thus cannot be both in $S$.

\msk

\nin
(6)\quad  
Observe that $\deg_{G^2}(v) \leq \deg(v) + \deg(v)(\Delta (G)-1) = \deg(v)\Delta(G)$ holds for every vertex $v$.
Hence, applying (5) and the Caro--Wei inequality, $\sid(G)  \geq   \alpha(G^2)  \geq   \sum_{v\in V} 1/(\deg_{G^2}(v) +1 ) \geq   \sum_{v\in V} 1/(\deg_G(v)\Delta(G)   +1) \geq  n/ (d(G)\Delta (G) +1)$, the last step obtained via Jensen's inequality. 

Alternatively,

 $$d(G^2)  = 2e(G^2)/n \leq  \bigl(\sum_{v\in V} \deg_{G^2}(v)\bigr)/n \leq  \Delta(G) \bigl(\sum_{v\in V} \deg(v)\bigr)/n = \Delta (G)d(G)$$
 and, by Tur\'an's theorem,

 $$\sid(G)  \geq   \alpha (G^2) \geq  n/(d(G^2)+1) \geq  n/ (d(G)\Delta (G) +1) \geq  n/( \Delta ^2(G)+1) \,.$$
\epf

\bpn
If $G$ is a bipartite graph with all degrees odd, then

\begin{itemize}
 \item[$(i)$] $\soc(G)=2$ and $\sid(G)\geq |G|/2\geq (\alpha(G)+1)/2$;
 \item[$(ii)$] if $G$ is also $d$-regular, with $d$ odd, then $\sid(G)=\alpha(G)=|G|/2$.
\end{itemize}
\epn

\bpf
By the degree assumption, each vertex class is an \sois, and one of them
 contains at least half of the vertices, hence proving $(i)$.
(The lower bound $(\alpha(G)+1)/2$ on $|G|/2$ is valid since $G$ cannot be edgeless.
Tightness is shown by the stars $K_{1,n-1}$, $n$ even.)

If $G$ is regular and bipartite, then it contains a perfect matching, consequently $\alpha(G)$ equals half the order of the graph.
This proves $(ii)$.
\epf

For $d$-regular graphs with $d$ even,  we have the following upper bound.

\bpn
\label{p:d-reg}
Assume that $d$ is even.
If $G$ is a $d$-regular graph, then $\sid(G) \leq \frac{d-1}{2d-1}\cdot|G|$,
 and this upper bound is tight for every $d\geq 2$
  in the class of $d$-regular connected graphs
 of any order $n$ of the form $n=t\cdot(2d-1)$.
Moreover, the bound is tight also in the class of $d$-regular bipartite graphs.
\epn

\bpf
Upper bound:
Let $A$ be an \sois\ of size $|A| = \sid(G)$, and let $B$ be the set of
 vertices having at least one neighbor in $A$.
Consider the bipartite graph $H$ formed by the $A$--$B$ edges.
Every $u\in A$ has degree $d$ in $H$ as $A$ is independent, and every
 $v\in B$ has degree at most $d-1$ in $H$, because it must be odd
 and $d$ is even.
Thus, $d\cdot |A| \leq (d-1)\cdot |B|$, and consequently
 $|G| \geq |A| + |B| \geq |A| + \frac{d}{d-1}\cdot|A|$, implying the assertion.

Lower bound:
Consider $t\geq 2$ copies of $K_{d,d-1}$, and take a perfect matching in the
 union of the $t$ $d$-element vertex classes in such a way that the
 obtained graph is connected.
(If $t=1$, just take a matching of $d/2$ edges inside the vertex class
 of size $d$.)
In particular, for $d=2$ this means $C_{3t}$.

To obtain bipartite examples, let $t=2$ and take the mentioned perfect
 matching between the two $d$-element vertex classes, keeping those two
 classes as independent sets.
\epf

We can also state an asymptotic version of the above assertion, for
 graphs that are not far from being regular.

\bpn
\label{p:regbip-asymp}
Let $\Delta$ be an even natural number, and for $j=1,2,\dots$ let
 $G_j$ be a graph on $n_j$ vertices, $n_j \to\infty$,
 of maximum degree $\Delta$ in which the number $z_j$ of vertices
 of degrees smaller than $\Delta$ satisfies $z_j = o(n_j)$.
Then $\limsup \sid(G_j) / n \leq \frac{\Delta-1}{2\Delta-1} \, n$.
\epn

\bpf
Let $B_j$ be a maximum odd independent set of $G_j$, $Z_j\subset V(G_j)$ the set of
 vertices of degree smaller than $\Delta$, and $B  =  B_j  \smin Z_j$.
Clearly, $|B_j|\geq n/(\Delta^2+1)$ and if $n_j$ is large enough then
 $B \neq \es $ as $z_j = o(n_j)$.
Due to the parity conditions, every vertex not in $B_j$ is adjacent to
 at most $\Delta-1$ vertices in $B$.
The number of edges emanating from $B_j$, say $e_j$, can be
 estimated as

  $$
    \Delta\cdot |B| = \Delta\cdot ( |B_j| - z_j ) \leq e_j \leq
     (\Delta-1)( n_j  - |B_j|) = (\Delta-1)( n_j - |B| - z_j) \,.
  $$
This gives $(2\Delta-1)|B| \leq (\Delta-1)(n_j - z_j )$ and
 $\sid(G_j) \leq |B| + z_j \leq \frac{\Delta-1}{2\Delta-1}\,(n_j-o(n_j))
 + o(n_j)$.
Hence $\limsup \sid(G_j) / n \leq \frac{\Delta-1}{2\Delta-1}\,n$.
\epf

Another extension, inspired by strongly regular graphs, can be stated as follows.

\bpn
Let $G$ be a $\Delta$-regular graph on $n$ vertices, in which every two
 adjacent vertices have at least $\lambda$ common neighbors.

\begin{itemize}
  \item[$(i)$] If $\Delta - \lambda \equiv 0 \ (\mathrm{mod} \ 2)$, then
   $\sid(G)  \leq \frac{\Delta - \lambda -1}{2\Delta - \lambda -1}\, n$.
  \item[$(ii)$] If $\Delta - \lambda = 1 \ (\mathrm{mod} \ 2)$, then
   $\sid(G) \leq \frac{\Delta - \lambda}{2\Delta - \lambda}\, n$.
\end{itemize}
\epn

\bpf
Let $B$ be a maximum odd independent set in $G$.

\msk

\nin
$(i)$\quad
Consider $v \in V\smin B$.
Either $v$ is not adjacent to any vertex
 in $B$, or it has at least one neighbor $u \in B$.
In the latter case none of the (at least) $\lambda$ common neighbors
 of $u$ and $v$ can be in $B$ (as they are neighbors of $u$).
So, $v$ can be adjacent to at most $\Delta - \lambda$ vertices in $B$,
 but in fact this means only $\Delta-\lambda-1$ because
 $\Delta-\lambda$ is even.
Hence $v$ is adjacent to at most $\Delta - \lambda -1$ vertices in $B$.

Now count the number of edges between $B$ and $V \smin B$ in two ways
 and find $\Delta \cdot |B| \leq (\Delta - \lambda -1)(n - |B|)$.
Hence $(2\Delta - \lambda -1)|B| \leq ( \Delta - \lambda -1 )n$,
  therefore $\sid(G) = |B| \leq \frac{\Delta - \lambda -1}{2\Delta - \lambda -1}\, n$,
 as claimed.

\msk

\nin
$(ii)$\quad   
The same proof applies, observing that in this case $v$ can be adjacent
 to at most $\Delta - \lambda$ members of $B$.
\epf

The particular case $\lambda = 0$ of $(ii)$ yields the upper bound
 $n/2$, which is tight e.g.\ in odd regular bipartite graphs.

Similar upper bounds concerning $K_{1,r}$-free
 graphs can be formulated.
We give them in detail in the companion part II of this paper  \cite{part-2} where it is used for various grids and chessboard  graphs, which are $K_{1,5}$-free.

We also observe that graphs for all feasible combinations of the
 parameters exist.

\bpn
Graphs $G$ of order $n$ with $\sid(G) = \alpha(G) = k$ exist under any of the conditions below:
  \begin{itemize}
    \item[$(i)$] $1 \leq k \leq n$;
    \item[$(ii)$]  $1 \leq k \leq n$, $k$ odd, $G$ connected;
    \item[$(iii)$] $k \mid n$ with $k$ odd, $G$ regular and connected;
    \item[$(iv)$] $2 \leq k \leq n-2$, $k$ even, $G$ connected.
  \end{itemize}
\epn

\bpf
For the above cases, we can take the following graphs.

\msk

\nin
$(i)$\quad
Take $G=K_{n-k+1} \cup (k-1) K_1$, or more generally any graph consisting of $k$ complete components.

\msk

\nin
$(ii)$\quad
The graph $G=K_{n-k} + E_k$, the split graph with complete join works;
and for $k<n$ also any complete $t$-partite graph $K_{a_1,\dots,a_t}$ where $a_1 \geq a_2 \geq \dots \geq a_t$, with $\sum a_j = n$ and $a_1=k$ odd.

\msk

\nin
$(iii)$\quad
Take $G=K_{k,k,\dots,k}$, the complete equipartite graph with $n/k$ partite classes of size $k$.

\msk

\nin
$(iv)$\quad
If $k=n-2$, take $K_2 + E_{n-2} - 2K_2$.
More generally, for any even $k\leq n-2$, begin with the complete join $K_{n-k} + E_{k-1}$ and supplement it with a further vertex $z$.
Delete an edge $vw$ where $v\in K_{n-k}$ and $w\in E_{k-1}$; and insert the new edge $vz$.
(For $n=4$, both ways yield $P_4$.)
\epf

Note that $k=n-1$ with $k$ even does not admit a connected graph with $\sid=k$, because the only connected graph having $\alpha = n-1$ is the star of degree $n-1$.
And, of course, $k=\alpha=n$ implies the union of singleton components, always disconnected.

\subsection{Forbidden pair lemma}

Here we state a property which is simple but makes many situations concerning odd independent sets more transparent.
We introduce two notions.

\bdf   \nev{Forbidden pair}
Two nonadjacent vertices $x,y$ form a \emph{forbidden pair} if they have a common neighbor $z$ such that $N[z]\subseteq N[x]\cup N[y]$.
\edf

\bdf   \nev{Forcing pair}
Two nonadjacent neighbors $x,y$ of a vertex $z$ form a \emph{forcing pair} if, for any vertex $w$ in $N(z)$ independent of $x$ and $y$, either $x,w$ or $y,w$ or both are forbidden pairs.
\edf

\blm
\label{l:force}
If $S\subseteq V(G)$ is an \sois, then it contains no forbidden pair, nor a forcing pair.
\elm

\bpf
Suppose for a contradiction that $x,y\in S$, where $x,y$ is a forbidden pair, say in regard to a common neighbor $z$, and $S\subseteq V(G)$ is odd independent.
By definition, $S$ should contain at least one further neighbor, say $w$, of $z$.
However, every neighbor of $z$ belongs to $N[x]\cup N[y]$, hence at least one of $wx\in E(G)$ and $wy\in E(G)$ holds, which is not allowed if $S$ is an independent set.

Assume next that $x,y\in S$ form a forcing pair, say in regard to a common neighbor $z$.
Due to odd independence, $S$ has to contain at least one further neighbor, say $w$, of $z$.
But then the forbidden pair $x,w$ or $y,w$ is contained in $S$, and we have just proved that this is not possible.
\epf

\btm
\label{t:clawfree}
If $G$ is claw-free, then $\sid(G) = \alpha(G^2)$ and $\soc(G) = \chi(G^2)$.
\etm

\bpf
We prove that a set $S\subseteq V(G)$ is odd independent in $G$
 if and only if $S$ is an independent set in $G^2$.
This will verify both equalities.

The ``if'' part is clear.
For ``only if'', consider any independent set $S$ in $G$, and
 suppose that $S$ is not independent in $G^2$.
Say, $x,y\in S$ where $xy\notin E(G)$  and $x,y\in N(v)$
 for some vertex $v$.
Since $G$ is claw-free, the pair $x,y$ dominates all neighbors of $v$.
Consequently $N[v]\subseteq N[x] \cup N[y]$.
Thus $x,y$ is a forbidden pair, and $S$ cannot be an \sois.
\epf

An important consequence is the following complexity result.

\btm
The problem of determining $\sid(G)$ is \NP-hard.
It also remains \NP-hard when the input graph is restricted
 (among others) to the line graphs $L(H)$ of graphs from
 any of the following graph classes:
 planar graphs of maximum degree $4$, bipartite graphs of
  maximum degree $3$, $C_4$-free bipartite graphs,
  $r$-regular graphs with $r \geq 5$, line graph graphs,
 chair-free graphs, or Hamiltonian graphs.
\etm

\bpf
We apply reduction from the problem Maximum Induced Matching (MIM),
 which asks, for a generic input graph $H$, to determine the largest
 number of edges no two of which are adjacent or incident in $H$.
Induced matching in $H$ precisely means a set of vertices in the
 line graph $L(H)$ which have pairwise distances at least 3,
 hence an independent set in the graph $(L(H))^2$.
Every line graph is claw-free, thus denoting $G:=L(H)$ we have
 $\sid(G)=\alpha(G^2)$ by Theorem \ref{t:clawfree}.
The number of steps in this reduction from $H$ to $G$ and $G^2$ is
 polynomial in terms of input size.
As a consequence, if MIM is intractable on a certain class of graphs
 $H$, then so is the determination of $\sid$ on the line graphs
 $L(H)$ of those graphs.
\NP-hardness of MIM on the graph classes listed in our theorem
 follows from the results proved in the references
 \cite{KS-94},  \cite{KR-Alg03}, \cite{L-IPL02}.
\epf

\brm
The overall presence of forbidden pairs characterizes claw-free graphs.
Namely, a graph is claw-free if and only if any two vertices at distance~$2$ form a forbidden pair for each of their common neighbors.
Indeed, if $G$ contains a claw with center $v$ and leaves $x,y,z$
 then $z\in N[v]\smin (N[x] \cup N[y])$,
 hence $x,y$ is not a forbidden pair.
\erm

Further applications of forbidden and forcing pairs are given in \cite{part-2}.

\subsection{Applications of $\soc(G)$ to lower bounds on $\sid(G)$}

In this short subsection we collect the main results from the papers \cite{CPST,GK-etal,KP-a24,P-a25}
 concerning upper  bounds on $\soc(G)$ and their implications to obtain lower bounds for $\sid(G)$ via the relation $\sid(G) \geq |G|/\soc(G)$.
This can also serve as a comparison to future results that will be obtained by direct research on $\sid(G)$.

For the next part (A) we use the standard notation
$\Mad(G) := \max \, \{ d(H) \mid H\subseteq G \}$ for the maximum average degree of $G$.

\msk

\nin
(A)\quad
\underline{From the paper \cite{KP-a24}}
\begin{itemize}
 \item Theorem 1.2: \ If $G$ is a graph with $\Mad(G) \leq  20/7$,
  then $\soc(G) \leq  \Delta (G)+4$.
 \
   $\Longrightarrow$
   \
     $\sid(G)  \geq  n/ ( \Delta (G) +4)$.
 \item Theorem 1.3: \ If $G$ is a graph with $\Delta (G) \geq  4$ and $\Mad(G) \leq  30/11$, then $\soc(G) \leq  \Delta (G)+3$.
 \
   $\Longrightarrow$
   \
     $\sid(G)  \geq  n/(\Delta (G) +3)$.
 \item Theorem 1.4: \ If $G$ is a $C_4$-free subcubic graph with $\Mad(G) \leq  30/11$, then $\soc(G) \leq  6$.
 \
   $\Longrightarrow$
   \
     $\sid(G)  \geq  n/ 6$.
 \item Corollary 1.5: \ For any planar graph $G$,
  $\soc(G) \leq  \Delta (G)+4$ if $g(G) \geq  7$, and
  $\soc(G) \leq  \Delta (G)+3$ if $g(G) \geq  8$.
 \
   $\Longrightarrow$
   \
     $\sid(G) \geq n/ (\Delta(G) +4)$ if $g(G) \geq  7$, and
     $\sid(G) \geq n/ (\Delta(G) +3)$ if $g(G) \geq  8$.
    These estimates are better than the one in our \cite{CPST}, $n/388$,
      for $\Delta \leq  384$ or $\leq 385$, respectively.
\end{itemize}

\msk

\nin
(B)\quad
\underline{From the paper \cite{GK-etal}}

\begin{itemize}
 \item Proposition 3: \ For every outerplanar graph $G$, $\soc(G) \leq  8$.
 \
   $\Longrightarrow$
   \
     $\sid(G)  \geq  n/ 8$.
 \item Theorem 4: \ For every proper minor-closed graph class $\cF$, there exists a constant $c(\cF)$ such that for every graph $G \in \cF$ we have $\soc(G) \leq  c( \cF)$.
   \
   $\Longrightarrow$
   \
     $\sid(G)  \geq  n/ c(\cF)$.
 \item Theorem 10: \ There exists a function $f$ such that for every
  graph $G$ we have $\soc(G) \leq  f(\tw(G))$, where $\tw(G)$ denotes
  the treewidth of $G$.
   \
   $\Longrightarrow$
   \
     $\sid(G)  \geq  n/ f(\tw(G))$.
 \item Theorem 19: \ There exists a function $f$ such that for every
  graph $G$ we have $\soc(G) \leq  f(\rtw(G))$, where $\rtw(G)$ denotes
  the row treewidth of $G$.
   \
   $\Longrightarrow$
   \
     $\sid(G)  \geq  n/ f(\rtw(G))$.
\end{itemize}
However, the implied constants in Theorems 4, 10, and 19 are not made explicit.

\msk

\newpage

\nin
(C)\quad
\underline{From the paper \cite{P-a25}}

\begin{itemize}
 \item Generalizing the main results from \cite{GK-etal}, the following major result is proved.
 
 Theorem 1.1: \ Every graph class of bounded expansion has bounded strong odd chromatic number;
i.e., for any graph class $\cF$ of bounded expansion there is a constant
$c(F)$ such that, for every graph $G \in \cF$, $\soc(G) \leq c(F)$.
   \
   $\Longrightarrow$
   \
     $\sid(G) \geq |G|/c(F)$.
\end{itemize}

\msk

\nin
(D)\quad
\underline{From our paper \cite{CPST}}

\begin{itemize}
 \item Proposition 2.1: \ For every tree $T$ we have $\soc(T) \leq  3$,
  and $\soc(T) = 2$ holds if and only if $T$ is odd (i.e., if all its
  vertices have odd degrees).
   \
   $\Longrightarrow$
   \
     $\sid(T)  \geq  n/ 3$  and $n/2$, respectively.
 \item Proposition 2.2: \ If $G$ is a connected unicyclic graph other than $C_5$, then $\soc(G) \leq  4$
   \
   $\Longrightarrow$
   \
     $\sid(G)  \geq  n/4$.
 \item Corollary 3.5: \ For every planar graph $G$, $\soc(G) \leq 388$.
   \
   $\Longrightarrow$
   \
     $\sid(G) \geq n/388$.
\end{itemize}

\section{Odd independence in some special graphs}
\label{s:spec-grs}

In this section we compute, via distinct approaches, the odd independence number for various classes of well-known graphs.

\subsection{Paths and cycles}

 The values $\soc(P_n)$ (equal to 3 for all $n\geq 3$) and $\soc(C_n)$
 (at most 4 unless $n=5$) are determined for all $n$
 in Section 2 of \cite{CPST}.
The corresponding values of \soin\ are given in the next assertion.

\bpn
\label{p:path-cyc}
We have
$\sid(P_n)=\lceil n/3\rceil$
 and
  $\sid(C_n)=\lceil(n-2)/3\rceil$.
\epn

\bpf
For paths we prove the recursion $\sid(P_n) = \sid(P_{n-3}) + 1$ for
 $n\geq 4$.
This, together with the obvious fact $\sid(P_1) = \sid(P_2) = \sid(P_3) = 1$,
 implies $\sid(P_n)=\lceil n/3\rceil$ by induction.

Let $P_n=v_1\cdots v_n$ be the path of order $n$, and let $S$ be an
 \sois\ of size $|S|=\sid(P_n)$ such that the sum
 $\sum_{v_i\in S} i$ is maximum.
Consider the vertex $v_k\in S$ of highest index in $S$.
It must be the case that $k=n$, for otherwise replacing $v_k$ with $v_n$
 we would obtain a set $S'$ such that $|S'|=|S|$ and also $S'$ is \soi,
  but the sum of indices in $S'$ is larger than that in $S$,
 contradicting the choice of $S$.
Hence $v_n\in S$ and $v_{n-1},v_{n-2}\notin S$.
Since $\{v_{n-1},v_{n-2}\}$ separates $v_n$ from the rest of the path
 in the strong sense that they have no common neighbors,
 $S$ is \soi\ in $P_n$ if and only if so is $S\smin\{v_n\}$
 in $P_{n-3}=v_1\cdots v_{n-3}$.
This implies the assertion.

For very short cycles it is clear that $\sid(C_3)=\sid(C_4)=\sid(C_5)=1$.
If $n>5$, consider $C_n=v_1\cdots v_n$, and let $S$ be a largest \sois.
Pick any $v_i\in S$.
Then $v_{i-1},v_{i-2},v_{i+1},v_{i+2}\notin S$ and $S$ is \soi\
  in $C_n$ if and only if so is $S\smin\{v_i\}$ in
 $C_n - \{v_{i-2},v_{i-1},v_i,v_{i+1},v_{i+2}\} \cong P_{n-5}$
 (subscript addition taken modulo $n$).
Thus, $\sid(C_n)=1+\sid(P_{n-5})$ holds, and the assertion follows by
 the formula on paths.
\epf

\subsection{Moore graphs}

Next we consider the famous $d$-regular graphs of diameter 2
 for which the Moore bound is tight.
It means $C_5$ for $d=2$, the Petersen graph for $d=3$, and the
 Hoffman--Singleton graph \cite{HS60} for $d=7$.
(The existence of a tight construction with $d=57$ on $57^2+1=3250$
 vertices is still open.)
Of course $\sid(C_5)=1$ and $\soc(C_5)=5$.
The other cases are more complicated.

We shall apply the following simple observation.

\blm
Assume that $G$ has girth at least $5$.
If $v\in V(G)$ is any vertex and $S\subseteq N(v)$ any set of odd size,
 then $S$ is an \sois\ in $G$.
\elm

\bpf
By assumption $v$ has an odd number $|S|$ of neighbors in $S$.
On the other hand, any $w\in V(G)\smin (S\cup\{v\})$ can have at most
 one neighbor in $S$, for otherwise it would create a $C_4$ together
 with $v$ and two vertices from $N(v)\cap N(w)$.
\epf

For simplicity, let us say that an \sois\ $S$ is trivial if $|S|\leq 1$,
 and non-trivial otherwise.

\bpn
\label{p:P}
The Petersen graph has $\sid=3$ and $\soc=6$.
Moreover, every non-trivial \sois\ is the neighborhood $N(v)$
 of a vertex $v$.
\epn

\bpf
Let $x,y\in S$ be two vertices in any non-trivial \sois\ $S$.
Since the diameter is 2, $x$ and $y$ have a common neighbor $w$.
Then, by parity reason, $S$ has to contain the third neighbor, say $z$,
 of $w$.
In this case, however, $N(x)\cup N(y)\cup N(z)$ is the entire vertex set,
 hence no more vertices can occur in $S$.
This implies $\sid=3$.

Consider now a \sokk\ with $\soc$ colors.
Say, $p$ vertex classes are trivial and $q$ are non-trivial.
We then have $p+q=\soc$ and $p+3q=10$, implying $q\leq 3$.
We observe that the case $q=3$ is impossible.
Indeed, $N(v)\cap N(v')=\es$ holds for two vertices $v,v'$ in the
 Petersen graph (or more generally
 in any graph of diameter 2) if and only if $v$ and $v'$ are adjacent.
Hence, $|S|=3$ would imply the presence of a triangle, while the
 Petersen graph has girth 5.
Consequently $q\leq 2$ and $p+q = 10-2q\geq 6$.

A \sokk\ with 6 colors can be achieved by taking the neighborhoods of
 two adjacent vertices and four trivial \sois s.
This completes the proof of $\soc=6$.
\epf

Let us recall that the Hoffman--Singleton graph \cite{HS60} is 7-regular, has diameter 2 and girth 5, and its order is $7^2+1=50$.

\bpn
\label{p:HS}
The Hoffman--Singleton graph has $\sid=15$ and $\soc\leq 20$.
On the other hand, it has no \sois\ with exactly $k$ vertices in the
 range $11\leq k\leq 14$.
\epn

\bpf
We label the vertices as $0,1,\dots,49$ in the way that $i\mapsto i+10$
 (mod 50) is an automorphism, and the first ten vertices have the
 following adjacencies:
\begin{center}
 \begin{tabular}{c|l}
  vertex & \,neighbors  \\
  \hline
  0 & \quad \ 1, \ 4, \ 13, \ 16, \ 26, \ 43, \ 49  \\
  1 & \quad \ 0, \ 2, \ \ \,6, \ 18, \ 28, \ 36, \ 47  \\
  2 & \quad \ 1, \ 3, \ \ \,8, \ 24, \ 34, \ 38, \ 45  \\
  3 & \quad \ 2, \ 4, \ 10, \ 19, \ 29, \ 40, \ 48  \\
  4 & \quad \ 0, \ 3, \ \ \,5, \ 22, \ 32, \ 35, \ 46  \\
  5 & \quad \ 4, \ 6, \ \ \,9, \ 12, \ 24, \ 27, \ 37  \\
  6 & \quad \ 1, \ 5, \ \ \,7, \ 14, \ 21, \ 30, \ 40  \\
  7 & \quad \ 6, \ 8, \ 11, \ 19, \ 25, \ 35, \ 49  \\
  8 & \quad \ 2, \ 7, \ \ \,9, \ 13, \ 22, \ 31, \ 41  \\
  9 & \quad \ 5, \ 8, \ 10, \ 17, \ 33, \ 43, \ 47
 \end{tabular}
\end{center}

We can generate an \sois\ of size 15 from the triplet $\{0,2,7\}$, taking
 $S=\{ 10j, \,10j+2, \,10j+7 \mid j=0,1,2,3,4 \}$.
Every vertex $i\notin S$ has exactly three neighbors\footnote{Those
 35 neighborhoods form a Kirkman triple system KTS$(15)$ over $S$.}
 in $S$.
Similar 15-element \sois s are obtained in the same way from the sets
 $\{0,5,8\}$, $\{1,3,5\}$, $\{1,4,9\}$, $\{2,6,9\}$ also.

Further, we can apply the idea of rotation by 10, that is applying $i\mapsto (i+10)$~mod~50, to construct a \sokk\
 with 20 colors.
Let $S_1=\{0,2,6,18,47\}$ and $S_2=\{1,3,24\}$.
These are \sois s, being subsets of odd cardinality in $N(1)$ and $N(2)$,
 respectively.
Rotating them with multiples of 10 we obtain 10 mutually disjoint \sois s
 which together cover 40 vertices.
The other 10 vertices can be taken as singleton color classes.

Suppose for a contradiction that there exists an \sois\ $S$ of size $k$,
 $11\leq k\leq 14$.
There are $7k$ edges with one end in $S$ and the other in its complement,
 which we denote by $Q$.
In this $S$--$Q$ bipartite graph all the $q:=50-k$ vertices $i\in Q$ have
 odd degrees.
The most balanced degree distribution of the ends of the $7k\geq 77$ edges among the $50-k\leq 39$ vertices of $Q$ can be determined in two steps.
First give 1 to each of the $50-k$.
There remain $8k-50$ to be distributed, but we must obtain odd degrees everywhere, so the increase cannot be smaller than 2.
Hence in the second step we give 2 to $\min\,\{50-k,4k-25\}=4k-25$ vertices of $Q$.
This yields $4k-25$ degrees 3 and $(50-k)-(4k-25)=75-5k$ degrees 1.
Each vertex $i\in Q$ having 3 has a neighborhood $N(i)\cap S$ that covers three vertex pairs in $S$; altogether this means at least $12k-75$.
Observe:
 \begin{center}
  \begin{tabular}{cc|cccccccccc}
   $k$ &  &  & 11 &  & 12 &  & 13 &  & 14 \\
   \hline
   $12k-75$ &  &  & 57 &  & 69 &  & 81 &  & 93 \\
   $\binom{k}{2}$ &  &  & 55 &  & 66 &  & 78 &  & 91 
  \end{tabular}
 \end{center}
Since $12k-75 > \binom{k}{2}$ holds in this range, there would occur vertex pairs in
 $S$ which are contained in the neighborhood of more than one vertex of $Q$.
This would imply that the graph contains $K_{2,2} \cong C_4$, a contradiction.
\epf

Also in the lower range, several cardinalities are impossible;
 for instance, the Hoffman--Singleton graph does not contain any
 2-element \sois\ $S$.
The reason is that if $u,v\in S$, then $u,v$ are not adjacent, hence they
 have a common neighbor, say $w$.
But then $w$ should have an odd number of further neighbors in $S$.

\subsection{Kneser graphs}

\btm
For every $k\geq 2$ and every $n\geq 2k$, the Kneser graph $\mathrm{KG}(n,k)$ has
 $\sid = \alpha$ if and only if $\binom{n-k-1}{k-1}$ is odd.
Equivalently, $\sid = \alpha$ if and only if
 $\binom{n-1}{k-1} \not\equiv \sum_{t=1}^k \binom{k}{t}\binom{n-k-1}{k-1-t}$ $(\mathrm{mod} \ 2)$.
\etm

\bpf
We first note that $\mathrm{KG}(2k,k)$ is a matching, i.e.\ regular of degree~1.
Thus, every independent set in $\mathrm{KG}(2k,k)$ is also odd independent and
 $\sid=\alpha$ holds.

Assume $n>2k$.
A largest independent set in $\mathrm{KG}(n,k)$ corresponds to a collection of the
 maximum number of pairwise intersecting $k$-element subsets of an
 $n$-element set.
According to the famous Erd\H os--Ko--Rado theorem, all such largest
 families $\cF$ are obtained by picking all $k$-sets containing a fixed element.
Every $k$-element set $F\notin \cF$ is disjoint from exactly
 $\binom{n-k-1}{k-1}$ members of $\cF$, and this is exactly the degree
 of vertices in $\mathrm{KG}(n,k)$.
This is odd if and only if the largest independent sets in $\mathrm{KG}(n,k)$ are
 also odd independent.

The second condition is equivalent to the first one because the equality $\sum_{t=1}^k \binom{k}{t}\binom{n-k-1}{k-1-t} = \binom{n-1}{k-1} - \binom{n-k-1}{k-1}$ is valid for all $n\geq 2k$.
\epf

\subsection{The graphs $S(K_n)$, $H_{n,n}$, $K_n \Boxx K_n$}

\nin
\underline{$S(K_n)$ --- the complete subdivision of $K_n$. }

\ssk

For $n\geq 2$ the graph $S(K_n)$ is obtained from $K_n$ by subdividing each edge with a new ``subdivision vertex''.
Hence it has $\binom{n}{2}+n = \binom{n+1}{2}$ vertices, $n$ of which have degree $n-1$ and $\binom{n}{2}$ have degree 2.
The $n$ vertices from $K_n$ will be referred to as ``original vertices''.

\bpn
If $n\geq 2$ is any integer, then:
  \begin{itemize}
      \item[$(i)$] $\soc(S(K_n))=n$ for all $n$, except that $\soc(S(K_2))=3$ and $\soc(S(K_4))=5$.
      \item[$(ii)$] $\sid(S(K_n))=\binom{n}{2}$ if $n$ is even, and $\sid(S(K_n))=\binom{n-1}{2}+1$ if $n$ is odd.
  \end{itemize}
\epn

\bpf
Denote the original vertices by $v_0,\dots,v_{n-1}$, and let $x_{i,j}$ be the subdivision vertex corresponding to the edge $v_iv_j$ for every $0\leq i<j\leq n-1$.
Observe that no \sois\ $S$ can contain more than one original vertex, because if $v_i,v_j\in S$ held, then $x_{i,j}$ would have exactly two neighbors in $S$.
This property immediately implies $\soc(S(K_n))\geq n$, and also $\sid(S(K_n))\leq \binom{n}{2}$ because a larger $S$ would require the simultaneous presence of all subdivision vertices together with one original vertex.

\msk

\nin
\emph{Proof of $(i)$ for odd $n$.}\quad
We label the color classes $S_i$ with $0,1,\dots,n-1\in \zzz_n$.
Define
$S_i = \{v_i\}\cup\{ x_{i+j,i-j} \mid j\in \zzz_n\smin \{0\} \}$, subscript addition taken modulo $n$.
Then the closed neighborhood of each original vertex contains each color from $0,1,\dots,n-1$ exactly once, and the two neighbors of each subdivision vertex have distinct colors, hence a \sokk\ with $n$ colors is obtained.

\msk

\nin
\emph{Proof of $(i)$ for even $n$.}\quad
We begin with noting that a \sokk\ with $n$ colors is not possible for $n=2$ and $n=4$.
Indeed, $S(K_2)\cong P_3$, which needs exactly three colors.
Moreover, assigning color $i$ to the original vertex $v_i$, we cannot use this color for any $x_{i,j}$ ($j\neq i$), so the possible choices of subdivision vertices of color $i$ correspond to a triangle in $K_4$, and at most one such appearance is possible (as any two triangles of $K_4$ share an edge).
Then the odd-degree condition with respect to $S$ implies that only one $x_{i,j}$ can get color $i$.
This means that at least two of the six subdivision vertices must have colors distinct from $0,1,2,3$.
On the other hand, this is feasible with five colors, e.g.\ assign color $i$ to $v_i$ and $x_{i+1,i+2}$ ($i=0,1,2,3$) and color 4 to $x_{0,2}$ and $x_{1,3}$.

Assume that $n\geq 6$.
We prove that $n$ colors suffice.
Also in this case we assign color $i$ ($i=0,1,\dots,n-1$) to the original vertex $v_i$.
The subdivision vertices are colored as follows:
  \begin{itemize}
      \item $x_{i,j}$ has color $n-1$ for all $0\leq i<j\leq n-3$.
      \item $x_{i,n-2}$ has color 0 for all $1\leq i\leq n-3$.
      \item $x_{i,n-1}$ has color $n-2$ for all $0\leq i\leq n-4$.
      \item $x_{0,n-2}$ and $x_{n-3,n-1}$ have color 1.
      \item $x_{n-2,n-1}$ has color $n-3$.
  \end{itemize}
This coloring corresponds to an edge decomposition of $K_n$ in which the edge classes are isomorphic to $K_{n-2}$, $K_{1,n-3}$, $K_{1,n-3}$, $2K_2$, and $K_2$, respectively.
Hence a \sokk\ of $S(K_n)$ is obtained.

\msk

\nin
\emph{Proof of $(ii)$.}\quad
If $n$ is even, then every original vertex is adjacent to an odd number of subdivision vertices, hence $\{x_{i,j} \mid 0\leq i<j\leq n-1\}$ is an \sois\ of size $\binom{n}{2}$.
If $n$ is odd, we take the \sois\ $\{v_0\}\cup\{x_{i,j} \mid 1\leq i<j\leq n-1\}$ of size $\binom{n-1}{2}+1$, and claim that no larger one exists.

Let $E(S)$ denote the set of edges in $K_n$ corresponding to the subdivision vertices belonging to $S$.
If $S$ is \soi, then all degrees in the graph induced by $E(S)$ are odd, hence this graph must have an even order.
This means at most $n-1$ because $n$ is odd.
Thus, $|E(S)|\leq \binom{n-1}{2}$.
In addition, $S$ can contain at most one original vertex, as told at the beginning of the proof.
This implies the validity of the assertion.
\epf

This construction yields that, contrary to 1-degenerate graphs, and also in sharp contrast with the classical chromatic number, $\soc$ is unbounded in the class of 2-degenerate graphs.

\bcr
There exist an infinite sequence of $2$-degenerate graphs $G$ such that $\soc(G)$ grows proportionally to $\sqrt{|G|}$.
\ecr

\bsk

\nin
\underline{$H_{n,n}$ --- the Half-graph.}

\ssk

The graphs $S(K_n)$ are examples of bipartite graphs with large $\sid$ and moderate $\soc$. 

The counterpart pathology is the half-graph $H_{n,n}$ with vertex set $A\cup B$, where $A = \{ u_1,\dots,u_n \}$, $B = \{ v_1,\dots,v_n \}$; vertex $u_i$ is adjacent to vertex $v_j$ for all $j \geq  i$. 

\bpn
$\soc(H_{n,n}) = n +1$ and $\sid(H_{n,n}) = 2$, all \sois s having
 the form $\{v_i , u_j \}$, $i<j$.
\epn

\bpf
Suppose for a contradiction that an \sois\ $S$ has more than one vertex
 in $B$.
Let $i < j$ be its two highest indices in $B$.
Then $u_i\in A$ is adjacent just to $v_i$ and $v_j$, contradicting the
 conditions on \sois s.
Hence $|S\cap B|=1$, and in the same way looking at the hypothetical
 two lowest indices in $A$ we also obtain $|S\cap A|=1$.
Consequently every \sois\ larger than a singleton is a nonadjacent
 vertex pair $v_i , u_j$, $i<j$.

Concerning any \sokk\ this implies that it is a color partition into an
 $A$--$B$ matching of non-edges, together with singletons.
Since $u_1$ and $v_n$ have degree $n$, at most $n-1$ non-singleton classes
 can occur, hence $\soc(H_{n,n}) > n$.
On the other hand, $n+1$ is achievable by taking $u_1$ and $v_n$ as
 singletons, with the further \sois s $\{u_i,v_{i-1}\}$ for $i=2,3,\dots,n$.
\epf

\bsk

\nin
\underline{$K_n \Boxx K_n$ --- Cartesian product of complete graphs.}

\ssk

This graph is opposed to $S(G)$, which has been seen to have large $\sid$
 and moderate $\soc$.

\bpn
\label{p: Kn-box-Kn}
For $H=K_n \Boxx K_n$,
$\sid(H) = 1$ and $\soc(H) = |H| \sim (\Delta(H))^2\!/4$.
\epn

\bpf
If $u,v\in H$ are nonadjacent vertices, then they have precisely two
 neighbors, say one of them is $w$.
However, $N(w)$ induces $2K_{n-1}$ in $H$, consequently $\{u,v\}$ cannot be
 extended to any independent set containing three vertices from $N(w)$.
This implies $\sid(H) = 1$.
Since $\sid(G) \soc(G) \geq |G|$ holds for every $G$, putting $G=H$
 and using the trivial inequality $\soc(H) \leq |H|$
 we also obtain $\soc(H) = |H|$.
The asymptotic formula follows as $H$ is $(2n-2)$-regular.
\epf

One can observe that the same proof gives $\sid(H) = 1$ and $\soc(H) = |H|$ more generally for $H=K_p \Boxx K_q$.
Note further that $G = K_p \Boxx K_q$ is claw-free and $G^2 \cong K_{pq}$, hence Proposition \ref{p: Kn-box-Kn} also follows directly from Theorem \ref{t:clawfree}.

\section{Automorphisms and recursive lower bound}
\label{s:recurs}

In this section we present a recursive method to obtain lower bounds
 on $\sid$ for graphs satisfying a certain symmetry property.
An application of this approach to cube graphs will be given.

For this purpose we define a binary operation $\mu$ acting on pairs of graphs.
Let $G$ be a vertex-labeled graph and $H$ any graph.
Define $\mu[G,H]$ as the graph obtained as follows.
Replace each vertex of $w\in V(H)$ with a copy $G_w=(V_w,E_w)$ of $G$,
 and between those copies which are adjacent in $H$,
 take a perfect matching that connects vertices of the same label.

\btm
\label{t:aut-bip}
Let $G$ and $H$ be graphs with the following properties:
 \begin{itemize}
  \item $G$ admits an automorphism $\eta$ such that $\eta(v)\in N(v)$
    for every $v\in V(G)$, and the length of each orbit of $\eta$ is even;
  \item $H$ is bipartite and each $w\in V(H)$ has an even degree.
 \end{itemize}
Then $\sid(\mu[G,H]) \geq |H|\cdot \sid(G)$.
\etm

\bpf
If $S$ is an \sois\ in $G$, then by automorphism so is $S^*:=\eta(S)$ as well.
Let $A,B$ be the two vertex classes of $H$, and consider an $S$ with
 $|S|=\sid(G)$.
If $w\in A$, let $S_w$ be the copy of $S$ in $G_w$; and if
 $w\in B$, let $S_w$ be the copy of $S^*$ in $G_w$.
We select $S^+ := \cup_{w\in V(H)} S_w$.
Clearly, $|S^+| = |H|\cdot \sid(G)$.
The proof will be done if we prove that $S^+$ is an \sois.

Trying to locate a hypothetical edge $u_1u_2$ inside $S^+$, there should
 exist $w_1,w_2$ such that $u_i\in S_{w_i}$ ($i=1,2$, $w_1\neq w_2$) and
 $w_1w_2\in E(H)$, because each $S_w$ is independent and there is no
 $G_{w_1}$--\,$G_{w_2}$ edge for $w_1w_2\notin E(H)$.
But then $u_1u_2$ should connect two copies of a $u\in V(G)$, which
 cannot be the case because one of them is in $A$ while the other is
 in $B$, hence $u_1 \mapsfrom u$ and $u_2 \mapsfrom \mu(u)\in N(u)$
 are nonadjacent in $\mu[G,H]$.

It remains to show that every vertex of $V(\mu[G,H])\smin S^+$
 adjacent to $S^+$ has an odd number of neighbors in $S^+$.
By symmetry reason, it suffices to consider vertices in the copy $G_{w_1}$
 of $G$.
Assume that $v_1\in (G_{w_1}) \smin S_{w_1}$ is the copy of $v\in V(G) \smin S$.

If $v$ has no neighbor in $S$, then referring to the inverse automorphism $\mu^{-1}$ we see that $\mu(v)$ has no neighbor in $S^*$, thus $v_1$ is not adjacent to $S^+$ either.
The other possibility is that $v$ has an odd number of neighbors in $S$.
This yields that $v_1$ has an odd number of neighbors in $S_{w_1}$.
All the other neighbors $u_i\in S^+$ of $v_1$, i.e.\ those in any $S_{w_i}$ where $w_i\in B$,
 are copies of the same vertex $u\in S^*$ where $u=\mu(v)$ holds by the definition of $\mu[G,H]$.
Since all degrees in $H$ are even, we obtain $|N(v_1) \cap S^+|
 \equiv |N(v_1) \cap S_{w_1}|$ (mod~2), hence odd.
It follows that $S^+$ is an \sois.
\epf

In a similar way we can also obtain:

\btm
\label{t:GxK2}
For any graph $G$, if $H$ is a bipartite graph with all vertex degrees
 odd, then
  $\sid(\mu[G,H]) \geq \frac{1}{2} |H|\cdot \sid(G\Boxx K_2)$.
\etm

\bpf
We apply notation and several ideas from the previous proof.
Let $A$ and $B$ be the two vertex classes of $H$, and $S\cup S^*$
 an \sois\ of maximum size $|S| + |S^*| = \sid(G\Boxx K_2)$
 in $G\Boxx K_2$, where $S$ is the subset in one copy of $G$ and $S^*$
 is the subset in the other copy.
We may assume $|A|\geq |B|$ and $|S|\geq |S^*|$ without loss of generality.

As in the previous proof, inside the copies $G_w \subset \mu[G,H]$ we
 select the copies of $S$ or $S^*$ if $w\in A$ or $w\in B$, respectively.
The argument in the previous proof clearly yields that the obtained set $S^+$
 is independent.
We need to prove that it is also odd independent.
Let $w_1w_2\in E(H)$ be any edge, $w_1\in A$, and let $v\in V(G)\smin S$
 be any vertex with at least one neighbor in $S\cup S^*$.
Consider the copy $v_1$ of $v$ in $G_{w_1}$.
The number of its neighbors in $V(G_{w_1})
 \cup V(G_{w_2})$ is odd, because the subgraph induced by this union
 is isomorphic to $G\Boxx K_2$ and $S\cup S^*$ is \soi.
The additional $\deg_H(w)-1$ neighbors of $w_1$ increase this with
$(\deg_H(w)-1)\cdot|N(v)\cap S^*|$, hence adding it to $|N(v)\cap (S\cup S^*)|$ the sum remains odd, due to the degree assumption on $H$.
The same is true for the copy $v_2$ of $v$ in $G_{w_2}$.
Thus, $S^+$ is an \sois.

We now have $\sid(\mu[G,H]) \geq |S^+| = |A|\cdot|S| + |B|\cdot |S^*| =
 \frac{1}{2}(|A|+|B|)|S| + \frac{1}{2}(|A|-|B|)|S| + |B|\cdot |S^*| \geq
 \frac{1}{2}(|A|+|B|)|S| + \frac{1}{2}(|A|-|B|)|S^*| + |B|\cdot |S^*| =
 \frac{1}{2}(|A|+|B|)(|S| + |S^*|) = \frac{1}{2} |H|\cdot \sid(G\Boxx K_2)$.
\epf

\brm
Although the two lower bounds above have been derived using very similar
 approach, the estimates are independent in the qualitative sense, because
 $\sid(G\Boxx K_2)$ is not expressible as a function of $\sid(G)$.
This fact is demonstrated by the following simple examples:
 \begin{itemize}
  \item $\sid(K_2)=1$ and $\sid(K_2\Boxx K_2) = 1$.
    ($K_2\Boxx K_2 \cong C_4$.)
  \item $\sid(C_4)=1$ and $\sid(C_4\Boxx K_2) = 4$.
    ($C_4\Boxx K_2 \cong Q_3$ is 3-regular bipartite.)
 \end{itemize}
\erm

We can offer a further approach based on the following notion.
Due to its generality there may occur interesting applications of it, although we have not found any so far.

\bdf
Let $G$ be a graph and $B$ an \sois\ in $G$.
We say that $B$ is odd-suitable if there is a partition $V(G)=V_1\cup V_2$
 such that the sets $B_1:=B\cap V_1$, $B_2:=B\cap V_2$ and
 $B_1^*:=N(B_1)\cap V_2$, $B_2^*:=N(B_2)\cap V_1$ satisfy the
 following conditions for $i=1,2$\,:
  \begin{itemize}
   \item $N(B_i^*)\cap V_i = B_i$,
   \item the $B_i$--$B_i^*$ edges induce an odd bipartite graph.
  \end{itemize}
\edf

In the above theorems concerning $\sid(\mu[G,H])$ the $B_i$--$B_i^*$
 subgraphs are just matchings.
In this way the notion of odd-suitable sets opens the road to further
 control on $\sid$ in certain graph classes built recursively.
Currently we do not have applications of this kind, therefore we leave it
 for future research to develop the theory in this direction.

\section{Application to hypercubes}
\label{s:grid}

Let $Q_d$ denote the $d$-dimensional hypercube, whose vertex set is $\{0,1\}^d$, and two
 vertices are adjacent if and only if they differ in exactly one coordinate.
Equivalently, one can take the subsets of $\{1,\dots,d\}$ as vertices,
 two subsets being adjacent if and only if their symmetric difference
 is a single element.
In this representation the mapping $\eta$ with $F \leftrightarrow
 F\cup\{d\}$ for all $F\subseteq \{1,\dots,d-1\}$ satisfies the
 conditions of Theorem \ref{t:aut-bip}, all orbits having length 2.
Further, one may also view $Q_d$ as $\mu[Q_{d-1},K_2] \cong Q_{d-1}\Box
 K_2$ or $Q_{d-2}\Boxx C_4 \cong \mu[Q_{d-2},C_4]$
  (cf.\ Theorem~\ref{t:GxK2} with 1-regular $K_2$ and
  Theorem \ref{t:aut-bip} with 2-regular $C_4$).

Before turning to \sois s, let us note that the strong odd chromatic
 number of hypercubes is easy to determine as follows.

\bpn
$\soc(Q_d)=2$ if $d$ is odd, and $\soc(Q_d)=4$ if $d$ is even.
\epn

\bpf
If $d$ is odd, then the unique proper vertex 2-coloring of $Q_d$ is also
 a \sokk, hence the first formula is clear.
If $d$ is even, then $Q_d$ decomposes to two copies of $Q_{d-1}$ with
 a perfect matching betrween them.
Hence a proper coloring with colors $1,2$ in the first copy and with
 colors $3,4$ in the second copy is a \sokk\ with four colors.
Fewer colors are not enough, because the odd chromatic number of $Q_d$
 for even $d$ is equal to 4, as proved in \cite[Theorem 3.2]{CPST},
 and $\soc(G)\geq \xod(G)$ holds for all graphs $G$.
\epf

From Proposition~\ref{p:d-reg} we obtain:

\bcr
\label{c:cube-UB}
If $d$ is odd, then $\sid(Q_d) = \alpha(Q_d) = 2^{d-1}$.
If $d$ is even, then
 $\sid(Q_d) \leq \bigl(1 - \frac{1}{2d-1}\bigr) \cdot 2^{d-1}$.
\ecr

We do not know, however, whether or not the upper bound for the case of even $d$ is tight.
We can prove some partial results.

\bpn
\label{p:cube-LB}
The following lower bounds and equalities are valid.
\

\begin{itemize}
 \item[$(i)$] $\sid(Q_{4k}) \geq 2 \sum_{i=1}^k \binom{4k-1}{2i-1}$
   for all $k\geq 1$.
 \item[$(ii)$] $\sid(Q_4) = 6$,  $\sid(Q_6) = 24$, $\sid(Q_8) = 112$.
 \item[$(iii)$] 
   $\sid(Q_{d}) \geq \frac{7}{16}\cdot 2^d$
   for all even $d\geq 8$.
\end{itemize}
\epn

\bpf
We represent the vertices of $Q_d$ as subsets of $\{1,\dots,d\}$.

\msk

\nin
$(i)$\quad
To obtain an odd independent set achieving the claimed size,
we select those subsets of size $1,3,\dots,2k-1$ which do not contain the
 element $4k$, and those of size $2k+1,2k+3,\dots,4k-1$ which contain $4k$.

If a $2k$-element set does not contain $4k$, it has $2k$ neighbors of size
 $2k-1$ (omitting any one of its elements),
  and a unique neighbor of size $2k+1$ (inserting $4k$ into it);
 together these are $2k+1$ selected neighbors.
If a $2k$-element set contains $4k$, then it has $2k-1$ neighbors of size
 $2k-1$ (omitting any one element other than $4k$),
  and $2k$ neighbors of size $2k+1$ (inserting any one element from its
   complement); together these mean $4k-1$ selected neighbors in $Q_{4k}$.
Small sets of even size $2s < 2k$ do not contain $4k$, hence they have
 $2s$ neighbors of size $2s-1$ and $4k-2s-1$ neighbors of size $2s+1$;
  together these are $4k-1$.
Large even sets of size $4k-2s$ are complementary to small ones, so they
 also have $4k-1$ selected neighbors.
Thus, an \sois\ is obtained, which has the claimed number of elements.

\msk

\nin
$(ii)$\quad
The lower bound 6 is the first case of $(i)$, and 24 is obtained by the
 $\mu[G,H]$ construction with $H=C_4$.
The lower bound 112 was obtained via computer aid.
Also tightness has been verified by a computer code.
The binary sequences describing the characteristic vectors of the 112 sets are as follows:

\msk

00000010, 00000100, 00000111, 00001000, 00001011, 00001101, 00001110, 

00010000, 00010011, 00010101, 00010110, 00011001, 00011010, 00011100, 

00100000, 00100011, 00100101, 00100110, 00101001, 00101100, 00101111, 

00110001, 00110010, 00110111, 00111000, 00111011, 00111101, 00111110, 

01000000, 01000011, 01000101, 01001001, 01001010, 01001100, 01001111, 

01010001, 01010010, 01010100, 01010111, 01011011, 01011101, 01011110, 

01100001, 01100010, 01100100, 01100111, 01101000, 01101011, 01101110, 

01110000, 01110101, 01110110, 01111001, 01111010, 01111100, 01111111, 

10000000, 10000011, 10000101, 10000110, 10001001, 10001010, 10001111, 

10010001, 10010100, 10010111, 10011000, 10011011, 10011101, 10011110, 

10100001, 10100010, 10100100, 10101000, 10101011, 10101101, 10101110, 

10110000, 10110011, 10110101, 10110110, 10111010, 10111100, 10111111, 

11000001, 11000010, 11000100, 11000111, 11001000, 11001101, 11001110, 

11010000, 11010011, 11010110, 11011001, 11011010, 11011100, 11011111, 

11100011, 11100101, 11100110, 11101001, 11101010, 11101100, 11101111, 

11110001, 11110010, 11110100, 11110111, 11111000, 11111011, 11111101.\break

A lower bound 104 for $d=8$ can be obtained by the following construction,
 that has a transparent structure.
We take $\{1,2,...,8\}$ as the underlying
 set whose $2^8$ subsets are the vertices of $Q_8$.
Select all 4-element subsets, take those 2-element subsets which form the
 4-regular bipartite Tur\'an graph (16 pairs), the 16 corresponding
 6-element sets whose complements form the Tur\'an graph, and finally the
  empty set and the entire set $\{1,2,...,8\}$.
These are $70+16+16+1+1=104$ vertices in $Q_8$.

A triplet has 5 extensions to a selected quadruplet, and contains
 either 2 or 0 selected pairs, depending on whether it intersects both
 $\{1,2,3,4\}$ and $\{5,6,7,8\}$ or not.
Hence, in $Q_8$ the number of selected neighbors of a triplet is 7 or 5.
A singleton has 4 extensions to a selected pair and contains the
 empty set, hence its number of selected neighbors is 5.
The degrees of 5-sets and 7-sets are analogous.
Thus, the selected sets generate an \sois\ in $Q_8$.

\msk

\nin
$(iii)$\quad
This estimate follows by Theorem~\ref{t:aut-bip} with $H=C_4$.
\epf

\section{Complements of triangle-free graphs}
\label{s:co-triangfree}

\def \diam {\mathrm{diam} }

Here we investigate graphs $G$ whose complement $\overline{G}$ is triangle-free.
We will analyze the possible combinations of the diameter of such a graph and its complement.

Let us begin with the following two complexity observations.

\bpn
If $\overline{G}$ is triangle-free, then
 $\soc(G)$ can be determined in polynomial time.
\epn

\bpf
In any proper coloring (not only in strong odd, but also in non-restricted) of $G$,
 every color class is either a singleton or a vertex-pair, because
 $\alpha(G)=2$.
Hence we can detect all odd-independent vertex pairs by brute force,
 their number is $O(n^2)$.
These pairs will be the edges of an auxiliary graph $H$.
This $H$ has matching number $m$ if and only if $G$ has
 $\soc(G) = n-m$, where $2m$ vertices are
 covered by the $m$ matching edges and the
 remaining $n-2m$ vertices form singleton color classes.
So, we can construct $H$ in polynomial time and then determine
 its matching number efficiently, applying some of the well-known algorithms.
\epf

\bpn
In the family $\cF_k$ of all graphs with $\alpha(G) \leq k$, determining $\sid(G)$ and $\alpha(G^2)$ can be done in polynomial time for every fixed $k$.
\epn

\bpf
This can be done by checking all $j$-vertex subsets, $j = 1,\dots,k$, \
(a)~for independence, (b) for odd independence, and (c) for independence in $G^2$.
\epf

The family $\cF_2$ is the family of the complements of triangle-free graphs, as a graph $G$ has $\alpha(G) \leq 2$ if and only if its complement $\overline{G}$ is triangle-free.
So, in $\cF_2$ we expect that some structural results (beyond the immediate complexity result) are also possible.

\bsk

We now turn to the analysis of $\diam(G)$ and $\diam(\overline{G})$, and their relation to odd independence.
Under the triangle-free assumption on $G$ we are going to prove:

\begin{itemize}
 \item If $\diam(G) \geq 4$, then $\sid(\overline{G}) = \alpha((\overline{G})^2) = 1$ and $\soc(G) = |G|$.
 \item If $\diam(G) = 3$ and $\diam(\overline{G}) = 2$, then $\sid(\overline{G}) = \alpha((\overline{G})^2) = 1$ and $\soc(\overline{G}) = |G|$.
 \item If $\diam(G) =3$ and $\diam(\overline{G}) = 3$, then $\alpha(\overline{G}) = \sid(\overline{G}) = \alpha((\overline{G})^2) = 2$.
 \item If $\diam(G) = 2$ and $\diam(\overline{G}) = 2$, then $\sid(\overline{G}) = \alpha((\overline{G})^2) = 1$ and $\soc(\overline{G}) = |G|$.
\end{itemize}

\begin{itemize}
\item[(1)] \quad
$\diam(G)\geq 4$, $G$ triangle-free
\end{itemize}

In this case it is well known that $\diam(\overline{G})\leq 2$;
 see, e.g., \cite{BM-book}.
However, then in $\overline{G}$, no \sois\ can have two vertices because the diameter is at most 2 and, as the independence number is at most 2, we conclude that $\sid(\overline{G}) = \alpha((\overline{G})^2) = 1$ and $\soc(\overline{G}) = |G|$.

Hence this graph class provides another example for $\sid(G) = 1$ or, equivalently, for $\soc(G) = |G|$.

A large class of further examples will be obtained below at the end of part (3).

\begin{itemize}
\item[(2)] \quad
$\diam(G)=3$, $\diam(\overline{G})=2$, $G$ triangle-free
\end{itemize}

Then $\sid(\overline{G}) = \alpha(\overline{G}^2) =1$ and $\soc(\overline{G}) = |G|$.

\begin{itemize}
\item[(3)] \quad
$\diam(G)=3$, $\diam(\overline{G})=3$, $G$ triangle-free
\end{itemize}

For such graphs, two vertices realizing the diameter of $\overline{G}$ form both a maximum independent set and also an odd independent set.
Further, $\alpha(\overline{G}) = \sid(\overline{G}) = \alpha(\overline{G}^2) = 2$.

Furthermore, if $\diam(G)=3$, then we know $\diam(\overline{G})\leq 3$.
Below we give constructions for both $\diam(\overline{G})=3$ and $\diam(\overline{G})=2$.

\paragraph{Construction of triangle-free graphs with diameter 3 whose complements have diameter 3.}
\

Such graphs exists, the smallest one is the self-complementary $G=P_4$ with $\sid(G) = \alpha(G) = 2$, cf.\ Problem \ref{pm:sid=alfa} $(i)$.
Also the graph $G = C_6=v_1v_2\dots v_6$ supplemented with the edges $(v_1, v_4)$ and $(v_2,v_5)$ is triangle-free of diameter 3.
Its complement consists of two triangles connected by an edge, hence having diameter 3 and $\sid(\overline{G}) = \alpha(\overline{G}) = 2$ (but $\alpha(G)=3$).

More generally, consider $H = K_{n, m} -e$ (one edge omitted), with $n \geq m \geq 2$.
This $H$ is triangle-free and has $\diam(H) = 3$, forced by the vertices of the removed edge.
Its complement $G=\overline{H}$ is $K_n \cup K_m$ supplemented with one edge connecting the two complete parts, hence having diameter 3. Further, $\sid(G) = \alpha(G) = 2$ (cf.\ again Problem \ref{pm:sid=alfa}).

This construction can be extended further,
 taking $K_{n, m}$ ($n \geq m \geq 2$) and deleting $tK_2$ where $1 \leq t \leq m -1$. This is a triangle-free graph (in fact bipartite) of diameter 3.
 The complement is the union of $K_n$ and $K_m$ connected by $tK_2$ between them.
Hence there are vertices $u$ in $K_n$
and $v$ in $K_m$ where $v$ is not connected directly by any edge to $K_n$ while $u$ is not connected directly by any edge to $K_m$. Then the distance between $u$ and $v$ is 3. We take $u$ and $v$ as an odd independent set, which is also obviously a maximum independent set, establishing $\sid=\alpha$.

\paragraph{Construction of triangle-free graphs with diameter 3 whose complements have diameter 2.}
\

In this case $\sid(\overline{G}) = 1$ and $\soc(\overline{G}) = |G|$ holds (cf.\ Problem \ref{pm:soc=n}).

Let $G$ be a triangle-free graph with $\diam(G)=2$, $|G| \geq 3$, and $\Delta(G) \leq |G| - 2$.
We construct $H = G \Boxx K_2$, i.e., two vertex-disjoint copies of $G$ supplemented with a perfect matching between the corresponding vertices.
This $H$ is a triangle-free graph with diameter 3.

We claim that the complement $\overline{H}$ of $H$ has diameter 2.
In this graph the two copies of $\overline{G}$ are
 joined by a complete bipartite graph from which just a matching has been removed.
So, inside each $\overline{G}$, any two vertices have a common neighbor
 in the opposite copy of $\overline{G}$, because $|G| \geq 3$.

Consider two vertices $u,v$ in distinct copies of $\overline{G}$.
If they are copies of distinct vertices of $G$, then they are directly adjacent in $\overline{H}$.
Otherwise, if they are copies of the same vertex say $w$,
 we know that $\deg_G(w)\leq |G|-2$, hence $u$ has a neighbor $z$ in its
 $\overline{G}$ and then $uzv$ is an $u$--$v$ path of length 2.
Thus, $\diam(\overline{H})=2$ holds, indeed.

Observe that this method can be generalized to constructions of triangle-free graphs of any diameter $t \geq 3$ whose complements have diameter 2.
Just start with a triangle-free graph $G$ having $\diam(G)=t$ and take $H = G \Boxx K_2$.
The above argument verifies $\diam(\overline{H})=2$, hence covering case (1) in a general way.

\begin{itemize}
\item[(4)] \quad
$\diam(G)=2$, $G$ triangle-free
\end{itemize}

The smallest such graph is $C_5$.
We will prove that under the assumptions of (4) either $\diam(\overline{G})=2$ or $\overline{G}$ is disconnected.
So, in the following we exclude finite diameter greater than 2.

\begin{itemize}
\item[(4.a)] \quad
$\diam(\overline{G})\geq 4$
\end{itemize}

If $\overline{G}$ is connected and has diameter at least 4, then there are vertices $v_1,v_2,v_3,v_4,v_5$ forming an induced path.
But then $\{v_1 ,v_3 ,v_5\}$ is an independent set in $\overline{G}$, contradicting the assumption that $G$ is triangle-free.

\begin{itemize}
\item[(4.b)] \quad
$\diam(\overline{G})=3$
\end{itemize}

Suppose that the induced path $xuvy$ realizes diameter 3
 in $\overline{G}$.
If $x$ has another neighbor $z$, then $z$ must be adjacent to $u$
 as otherwise $\{z, u , y\}$ is an independent set of
 cardinality 3 in $\overline{G}$, which is forbidden.
Similarly, if $z$ is a neighbor of $y$, then it is also
 adjacent to $v$.
For the same reason both $N[x]$ and $N[y]$ induce complete subgraphs in $\overline{G}$.
What is more, the union of these two cliques is the entire vertex set,
 for otherwise a vertex $z$ not in their union would form a
 3-element independent set with $x$ and $y$.
Thus, $\overline{G}\supseteq (K_a\cup K_b) + e$ where $e$ is the edge $uv$
 and $a+b=|G|$ (namely, $a=\deg(x)+1$, $b=\deg(y)+1$).
But then $G\subseteq K_{a,b}-e$ holds,
 which implies $\diam(G)\geq 3$,
 hence contradicting the basic assumption of (4).

\paragraph{Construction of triangle-free graphs with diameter 2 whose complements have diameter 2.}
\

Start with $K_{n,m}$, $n \geq m \geq 2$, fix a matching of $t$ edges, $1 \leq t \leq m$, and subdivide these $t$ matching-edges into paths of length 2; call this graph $G = G(n,m,t)$.
This is a triangle-free graph of diameter 2.
(Note that $m > 1$ is mandatory, otherwise the graph would have diameter 3.)

The complement $H=H(n,m,t)$ of $G(n,m,t)$ is connected,
 thanks to $m>1$ and to the middle vertex of any subdivided edge.
Hence, by the analysis above, we know without further checking
 that $\diam(\overline{G})=2$.

\section{Graphs with $\soc(G)  = \Delta^2(G) +1$ and $\sid(G)\geq |G|/(\Delta^2(G) +1)$}
\label{s:maxdeg}

For a simpler notation, we will write $\Delta^2(G)$ for
 $(\Delta(G))^2$ throughout.

We have already seen that $\soc(G) \leq \chi(G^2) \leq \Delta(G^2) +1 \leq \Delta^2(G) +1$ holds for every graph $G$. 
So, $\chi(G^2) =\Delta^2(G) +1$ is a necessary condition for
 $\soc(G) = \Delta^2(G) +1$. 
This problem has a rich history; we refer to \cite{C-23} for the
 information needed here.

It is well known that a connected graph $G$ with $\Delta=\Delta(G)>2$ has
 $\chi(G^2) = \Delta^2(G) +1$ if and only if $G^2 = K_{\Delta^2(G) +1}$. 
As discovered by Hoffman and Singleton \cite{HS60}, such graphs exist
 if and only if $\Delta \in \{2 ,3 ,7\}$ and possibly $\Delta = 57$
 (whose existence is a major open problem in graph theory). 
We consider several cases.

\msk

(1)\quad If $\Delta  = 1$, then clearly $G=K_2$ and
 $\soc(G)=2=\Delta^2(G) + 1$.
We immediately infer that for $H=tK_2$,
 $\sid(H) = t = |H|/\soc(H) = |H|/(\Delta^2(H) + 1)$.

\msk

(2)\quad For $\Delta  = 2$ we already proved in \cite{CPST} that $\soc(G) \leq  5$, with equality if and only if at least one component is $C_5$.
Hence $\sid (G)\geq  |G|/5 = |G|/(\Delta ^2(G) +1)$, with equality if and only if  $G \cong tC_5$.  

\msk

(3)\quad $\Delta   = 3$.

It is proved in \cite{CR-16} that if $\Delta(G) \geq 3$ and $G$ is
 a connected graph other than the Petersen graph, the Hoffman--Singleton
  graph, or a graph with $\Delta(G) = 57$ and $G^2 = K_{57^2 +1}$,
    then $\chi(G^2) \leq \Delta^2(G) - 1$.
However, for $\soc(G) = \Delta^2(G) +1$, the Petersen graph is
 excluded by Proposition \ref{p:P},  and the
 Hoffman--Singleton graph is excluded by Proposition \ref{p:HS}.
Hence it remains to consider the case $\Delta = 57$, if it exists at all.   
To exclude this final possibility we need the following lemma.

\blm
Denote $\epsilon (\Delta) = 0$ if $\Delta$ is odd and $\epsilon (\Delta ) = 1$ if $\Delta$ is even.
If $G$ has girth at least $5$ and $\Delta(G)=\Delta$, then:
 \begin{itemize}
   \item[$(i)$] $\sid (G)  \geq   \Delta  - \epsilon (\Delta )$.
   \item[$(ii)$] $\soc(G)  \leq  n - ( \Delta  - \epsilon (\Delta )) +1$.
 \end{itemize}
\elm

\bpf
Assume that $G$ has girth at least $5$.
Let $v$ be a vertex with $\deg(v) = \Delta (G)$.
Consider $A = N(v)$, or one fewer vertex if $\Delta$ is even.
Note that $A$ is independent, since $G$ has no 3-cycles.
Vertex $v$ is adjacent to the entire $A$, whose size is odd.
Any other vertex can have at most one neighbor in $A$, because $G$ has
 no 4-cycles.
Hence $A$ is an odd independent set of size
 $\Delta - \epsilon(\Delta )$, proving~$(i)$.

Coloring all vertices of $A$ with the same color and all other
 vertices with distinct colors, we obtain $\soc \leq  n - ( 
\Delta  - \epsilon(\Delta)) +1$, as stated in $(ii)$. 
\epf	

A Moore graph $G$ with $\Delta = 57$, if it exists, is known to have girth 5.
Hence, by the above lemma, $\soc(G) \leq \Delta^2 - \Delta +2$, and so it is also excluded.  
We conclude:

\btm
\label{t:D2G}
If $\Delta(G) \geq 3$, then $\sid(G) \geq |G|/((\Delta(G))^2-1)$.
\etm

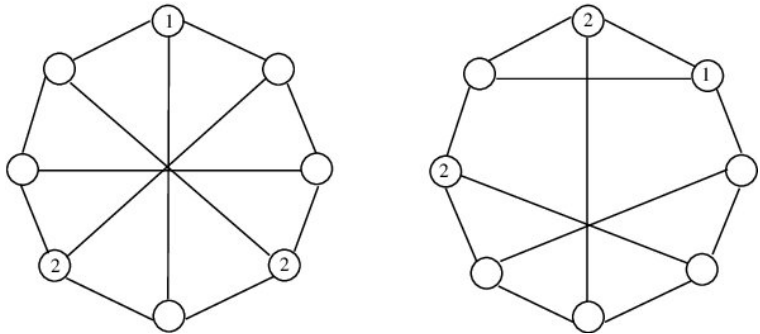
\begin{figure}[ht!]
\centering
\setlength{\tabcolsep}{1em}
\begin{tabular}{cc}

\begin{tikzpicture}[scale=1.2,
  every node/.style={circle, fill=black, draw=black, inner sep=2pt}]
\def\n{8}
\def\r{2}

\foreach \i in {1,...,\n} {
  \node (v\i) at ({90-360*(\i-1)/\n}:\r) {};
}

\foreach \i [evaluate=\i as \j using {int(mod(\i,\n)+1)}] in {1,...,\n} {
  \draw (v\i) -- (v\j);
}

\foreach \i in {1,...,4} {
  \pgfmathtruncatemacro{\k}{int(mod(\i+3,8)+1)}
  \draw (v\i) -- (v\k);
}
\end{tikzpicture}
&
\begin{tikzpicture}[scale=1.2,
  every node/.style={circle, fill=black, draw=black, inner sep=2pt}]
\def\n{8}
\def\r{2}

\foreach \i in {1,...,\n} {
  \node (v\i) at ({90-360*(\i-1)/\n}:\r) {};
}

\foreach \i [evaluate=\i as \j using {int(mod(\i,\n)+1)}] in {1,...,\n} {
  \draw (v\i) -- (v\j);
}

\draw (v1) -- (v5);
\draw (v2) -- (v8);
\draw (v3) -- (v6);
\draw (v4) -- (v7);
\end{tikzpicture}
\end{tabular}
\caption{Two 3-regular graphs of order 8 \label{fig:3-reg-order-8}}
\end{figure}

\brm
There exist some graphs with $\chi(G^2) = \Delta^2(G) - 1$.  
For $\Delta = 3$ this is the critical graph for the Ramsey number $R(4 ,3)$ on 8 vertices, for which one can check that 
$\soc(G) = 8$ and $\sid(G) = 1$.
Also a further graph with these parameters exists, see Figure \ref{fig:3-reg-order-8}.
Hence the bound in the above theorem is sharp for $\Delta = 3$. 
There exist also such graphs with $\Delta = 4$ and $\Delta = 5$, mentioned in \cite{C-23}; but for $\Delta \geq 6$ no example of such graphs is known so far.
\erm

\newpage

\section{Concluding remarks and open problems}
\label{s:concl}

In this section we offer a list of open problems for future research.

We know $\alpha(G) \geq \sid(G) \geq \alpha(G^2)$. Trivially for a connected graph $G$, $\alpha(G) = \alpha(G^2)$ if and only if $G = K_n$, $n \geq 1$.

\bpm
\label{pm:sid=alfa}
\
\begin{itemize}
\item[$(i)$]
Characterize the graphs $G$ for which $\sid(G) = \alpha(G)$ holds;
 or determine the complexity of deciding if $\sid(G) = \alpha(G)$.

\item[$(ii)$]
Characterize the graphs for which $\sid(G) = \alpha(G^2)$ holds;
 or determine the complexity of deciding if $\sid(G) = \alpha(G^2)$.
\end{itemize}
\epm

The first part has been done for families like the Kneser graphs, the Cartesian products
 $K_p \Boxx K_q$, and trivially for odd regular bipartite graphs.
The next assertion implies, however, that there does not exist an
 ``induced subgraph characterization'' for $\sid(G) = \alpha(G)$.

\bpn
For every graph $G$ there exists a graph $H$ such that $|H|=|G|+1$,
 $\sid(H) = \alpha(H) = \alpha(G) + 1$, and $G$ is an induced
 subgraph of $H$.
\epn

\bpf
Let $B$ be a largest independent set in $G$.
Construct $H$ by taking a new vertex $w$ and joining it to all those
 vertices of $V(G)\smin B$ which have an even number of neighbors in $B$.
(In particular, if $B$ is odd independent in $G$, then $H=G\cup K_1$.)
By construction it follows that $B^+:=B\cup \{w\}$ is an \sois, because if
 $v$ has an odd number of neighbors in $B$, then its degree
 remains unchanged, and otherwise it is increased by 1 due to the edge $vw$.
Hence
 $|B^+|\leq \sid(H)\leq \alpha(H)\leq \alpha(G)+1 = |B|+1 = |B^+|$.
Thus, equality holds throughout.
\epf

Note that the above construction of $H$ for a given $G$ needs to know a
 largest independent set in $G$ explicitly, which is algorithmically
 hard to find in general.

\bpm
\label{pm:soc=n}
Characterize graphs $G$ with $\soc(G) = |G|$, equivalently
 with $\sid (G) = 1$; or determine the complexity of the
  corresponding decision problem.
\epm

Clearly, a necessary condition is that $G$ has diameter at  most 2.
One can also make the following further related observations:
 \begin{itemize}
  \item[(a)] If $\soc(G) = |G|$, then $\soc(G +K _r) = |G| + r$.
  \item[(b)] If $G$ is claw-free, we know $\soc(G) = \chi(G^2)$ and  $\sid(G) = \alpha(G^2)$.
   Moreover, for any graph $G$ obviously $\chi(G^2) = |G|$ holds
    if and only if $G$ has diameter at most $2$.
   Hence, for claw-free graphs $G$, $\sid(G) = 1$ and $\soc(G) = |G|$ if and only if $G$ has diameter at most 2.
 \end{itemize}

\bpm
\label{pm:sid=reg}
Let $r$ be odd.
Characterize the connected triangle-free $r$-regular graphs for which
 \begin{itemize}
  \item[$(i)$] 
 $N(v)$ is a maximum \sois\ for every vertex $v$.
  \item[$(ii)$] 
 $N(v)$ is a maximum \sois\ for some vertex $v$.
  \item[$(iii)$] 
 $\sid (G) = r$.
 \end{itemize}
\epm

There are two known examples of such graphs, the  Petersen graph 
 and $K_{t ,t}$ where $t$ is odd.
Are there more?
We note that the two mentioned examples belong to the class of strongly
 regular graphs srg$(n,r, \lambda=0 , \mu =$\,odd).
The other five known members of this class, however, do not satisfy
 the requirements of Problem \ref{pm:sid=reg}.
(It is a famous open problem whether srg$(n,r, \lambda=0 , \mu)$
 contains further members beyond the seven known ones: namely, the
 pentagon, Petersen, Clebsch, Hoffman--Singleton, Gewirtz, Mesner-M22,
  and Higman--Sims.)

Observe that a necessary condition for $(i)$ and $(ii)$ is to have
 diameter at most 3.
Concerning $(iii)$ with $r \geq 3$ (not necessarily odd), since there is
 no request on any $N(v)$, we have $r = \sid(G) \geq |G|/(r^2 - 1)$,
 hence $|G| \leq r(r^2-1)$.

\bpm
Let $G$ be any outerplanar graph.
It is known \cite{GK-etal} that $\soc(G)  \leq  8$, moreover $\soc (G) =  7$ holds
 for $G = P_6 +K_1$.
Is it true that $\sid (G)  \geq  n/7$ for any outerplanar $G$\,?
\epm

The lower bound $\sid (G)  \geq n/7$ is valid at least for outerplanar graphs $G$ with $\Delta(G)\leq 6$, because then also $\chi(G^2)\leq 7$ holds, as proved in \cite[Corollary 3.13]{AH10}.

\bpm
In view of the lower bounds
 $\sid (G)  \geq  |G|/(\Delta^2(G) +1)$ and
 $\sid (G)  \geq  |G|/(\Delta^2(G) -1)$ for $\Delta \geq  3$, it is of
 interest to consider $c(\Delta ) := \min \{ c :$ there exists a graph
 such that $\sid (G) = c|G|/(\Delta^2(G) +1 ) \}$.
\epm

Here $c(1) = c(2) =  1$, but $c(3)  = 5/4$ realized by two graphs of
 order~8; one of them is the critical graph for $R(3,4)$.
The graphs $K_p \Boxx K_p$ and $K_p \Boxx K_{p+1}$ show that we cannot
 expect more than $c(\Delta )  = (4+o(1) )\cdot \Delta $.

To finish, for hypercubes the following problem remains open.

\bpm
Determine $\sid(Q_d)$ when $d$ is even, or improve upon the results in Corollary \ref{c:cube-UB} and Proposition \ref{p:cube-LB}.

In particular, is it true that if $d$ is even, then
 $\lim_{d\to\infty} \sid(Q_d)/2^d = 1/2$\,?
\epm

\paragraph{Acknowledgement.}

This research was supported in part
 by the Slovenian Research Agency ARIS program P1-0383, project J1-4351, and the annual work program of Rudolfovo (R. \v S.) and
 by the ERC Advanced Grant \mbox{``ERMiD''} (Zs.~T.).

\end{document}